\pdfoutput=1

\documentclass[a4paper, oneside]{article}
\usepackage{titlesec}
\titlespacing*{\section}
{0pt}{5.5ex plus 1ex minus .2ex}{4.3ex plus .2ex}

\titlespacing*{\subsection}
{0pt}{5.5ex plus 1ex minus .2ex}{4.3ex plus .2ex}

\usepackage[left=4cm,right=4cm,top=2.5cm,bottom=2.5cm,includeheadfoot]{geometry}

\usepackage{authblk}
\usepackage{makeidx}         
\usepackage{graphicx}   
\usepackage{multicol}        
\usepackage[bottom]{footmisc}

\usepackage{newtxtext}       %
\usepackage{newtxmath}       

\usepackage{caption}
\usepackage{units} 
\usepackage{tikz}
\tikzset{every path/.append style={line width=1pt}}
\usepackage{siunitx}
\usepackage{pgfplots}

\newtheorem{Algorithm}{Algorithm}

\newcommand{\keywords}{\textbf{Key words.}  }
\newtheorem{problem}{Problem}
\newtheorem{remark}{Remark}
\newcommand{\acknowledegment}{\textbf{Acknowledgments.}  }

\title{Adaptive Numerical Simulation of a Phase-field Fracture Model in Mixed Form tested on an L-shaped Specimen with High Poisson Ratios}
\author[1]{Katrin Mang}
\author[2]{Mirjam Walloth}
\author[1]{Thomas Wick}
\author[2]{Winnifried Wollner}
\affil[1]{Institute of Applied Mathematics, Leibniz Universit{\"a}t Hannover, Welfengarten 1, 30167 Hannover, Germany\\ 

katrin.mang@ifam.uni-hannover.de, thomas.wick@ifam.uni-hannover.de}

\affil[2]{Department of Mathematics, Technische Universit{\"a}t Darmstadt, Dolivostrasse 15, 64293 Darmstadt, Germany\\

walloth@mathematik.tu-darmstadt.de, wollner@mathematik.tu-darmstadt.de}
\date{}                     
\begin{document}

%
%
\maketitle

\abstract{
This work presents a new adaptive approach for the numerical simulation of a phase-field model for fractures in nearly incompressible solids. 
In order to cope with locking effects, we use a recently proposed mixed form where we have a hydro-static pressure as additional unknown besides the displacement field and the phase-field variable.
To fulfill the fracture irreversibility constraint, we consider a formulation as a variational inequality in the phase-field variable. 
For adaptive mesh refinement, we use a recently developed residual-type a posteriori error estimator for the phase-field variational inequality which is efficient and reliable, and robust with respect to the phase-field regularization parameter. 
The proposed model and the adaptive error-based refinement strategy
are demonstrated by means of numerical tests derived from the L-shaped panel test, originally developed for concrete.
Here, the Poisson's ratio is changed from the standard setting towards the incompressible limit $\nu \to 0.5$.
}



\section{Introduction}\label{sec_intro}
Crack propagation is one of the major research topics in mechanical, energy, and environmental engineering. A well-established variational approach for Griffith's~\cite{griffith1920phenomena} 
quasi-static brittle fracture was introduced by Francfort and Marigo~\cite{FraMar98}. 
Miehe et al.~\cite{MieWelHof10a} introduced the name `phase-field modeling' for this
variational approach.
If the observed solid is assumed to be nearly incompressible, the classical phase-field fracture model fails due to volume-locking.
In this work, we combine the mixed problem formulation, recently proposed by the authors in~\cite{mang2019phase}, with the adaptive numerical solution based on a residual-type error estimator for the arising phase-field variational inequality~\cite{mang2019mesh,Walloth:2018b}. 
This allows to simulate crack propagation on adaptive refined meshes in nearly incompressible materials by using the phase-field method.


\section{A Phase-field Model for Nearly Incompressible Solids}\label{NewModel}

\subsection{Notation and Spaces}
We emanate from a two-dimensional, open and smooth domain $\Omega\subset \mathbb{R}^2$. 
Let $I$ be a loading interval $[0,T]$, where $T>0$ is the end time value.
A displacement function $u:(\Omega \times I) \to \mathbb{R}^2$ is defined on the domain $\Omega$.
On a subset $\Gamma_D\subset\partial\Omega$ of the boundary, we enforce Dirichlet boundary conditions. 
For the phase-field variable $\varphi:(\Omega\times I)\rightarrow [0,1]$ with $\varphi=0$ in the crack and $\varphi=1$ in the unbroken material, we have homogeneous Neumann values $\nabla\varphi\cdot n= 0$ on the whole boundary, where $n$ is the unit outward normal to the boundary.
The physics of the underlying problem ask to enforce crack irreversibility, i.e., that $\varphi$ is monotone non-increasing with respect
to $t\in I$.
 
By $(a,b) := \int_\Omega a \cdot b\ dx$ for vectors $a, b$, the $L^2$ scalar-product is denoted. 
The Frobenius scalar product of two matrices of the same dimension is defined as $A:B:= \sum_i \sum_j a_{ij} b_{ij}$ and therewith the $L^2$-scalar product is given by $(A,B) := \int_\Omega A : B\ dx$. 

For a weak problem formulation, we consider a subdivision $0 =t_0 <\ldots < t_N = T$ of the interval
$I$. In each time step, we define approximations $(u^n,\varphi^n) \approx (u(t_n),\varphi(t_n))$ and hence the irreversibility condition is approximated by
$\varphi^{n}\leq\varphi^{n-1}$ for all $n = 1, \ldots, N$.
To simplify the notation, we omit the superscript $(\cdot)^n$ and set $u:=u^n$ and $\varphi:=\varphi^n$, whenever it is clear from the context. 
The phase-field space is $\mathcal{W} := H^1(\Omega)$ with a feasible set
$\mathcal{K}:=\{\psi\in \mathcal{W}\mid \psi\leq \varphi^{n-1}\leq 1\}$.
Further, we define the function spaces $\mathcal{V}:= (H_0^1(\Omega))^2:= \{w \in (H^1(\Omega))^2\mid w=0\mbox{ a.e. on }\Gamma_D\}$, $\mathcal{U}:=L_2(\Omega)$, and
$\mathcal{X}:=\{\Lambda\in \mathcal{W}^{\ast}\mid \Lambda\ge 0\}$, where $\mathcal{W}^{\ast}$ is the dual space of $\mathcal{W}$. Further, let $u_D\in (H^1(\Omega))^2\cap C^0(\Gamma_D)$ be a continuation of the Dirichlet-data.
For the classical phase-field fracture model, we refer to Miehe et al.~\cite{MieWelHof10a}. In the next section, the mixed form of the phase-field fracture model is formulated.

\subsection{Mixed Phase-field Fracture Model}\label{mixedModel}
The stress tensor $\sigma(u)$ is given by $\sigma(u) := 2 \mu E_{\text{lin}}(u) + \lambda \text{tr} (E_{\text{lin}}(u)) \textbf{I}$ with the Lam\'e coefficients $\mu,\lambda > 0$.
The linearized strain tensor therein is defined as $E_{\text{lin}}(u):=\frac{1}{2} (\nabla u + \nabla u^T)$. 
By $\textbf{I}$, the two-dimensional identity matrix is denoted. For a mixed formulation of the problem, we define
\begin{align*}
 p:= \lambda \nabla \cdot u,
\end{align*}
with $p \in \mathcal{U}$,
such that the pure elasticity equation reads as follows:\\
Find $u \in \mathcal{V}$ and $p \in \mathcal{U}$ such that 
 \begin{align*}
 \begin{aligned}
  2\mu (E_{\text{lin}}(u),E_{\text{lin}}(w))+(\nabla \cdot w, p)=&\ 0\quad \forall w \in\ \mathcal{V},\\
 (\nabla \cdot u,q) - \nicefrac{1}{\lambda} (p,q)=&\ 0\quad \forall q\in\ \mathcal{U}.
 \end{aligned}
 \end{align*}

Following~\cite{miehe2010phase}, we consider the tensile ($\sigma^+(u,p)$) and compressive ($\sigma^-(u,p)$) parts of the stress tensor. For this reason, the positive part of the pressure $p^+ \in L_2(\Omega)$ has to be defined as
$p^+ := \max \{p,0\}$, and $E_{\text{lin}}^+(u)$ is given as the 
projection of $E_{\text{lin}}(u)$ onto positive semidefinite matrices.
Now, we can split the stress tensor $\sigma(u,p)$ as:
\begin{align*}
\begin{aligned}
 \sigma^+(u,p):=&\ 2\mu E_{\text{lin}}^+(u) + p^+ \textbf{I},\\
 \sigma^-(u,p):=&\ 2\mu (E_{\text{lin}}(u)-E_{\text{lin}}^+(u))+ (p - p^+) \textbf{I}.
 \end{aligned}
\end{align*}

In the following, the critical energy release rate is denoted by $G_c$ and a degradation function is defined as $g(\varphi):=(1-\kappa)\varphi^2 + \kappa,$
with a small regularization parameter $\kappa > 0$. 
Next, we can formulate the mixed phase-field problem in incremental form~\cite{mang2019phase}:

\begin{problem}[Mixed Phase-field Formulation]\label{formMixed}
Given the initial data $\varphi^{n-1} \in \mathcal{K}$,
find $u:= u^n \in \{u_D+\mathcal{V}\}$, $p:=p^n \in \mathcal{U}$ and $\varphi:=
\varphi^n \in \mathcal{K}$ for loading steps
$n=1,2,\ldots, N$ such that 
 \begin{align*}
 \begin{aligned}
 g(\varphi^{n-1})(\sigma^+(u,p),E_{\text{lin}}(w) ) + (\sigma^-(u,p),E_{\text{lin}}(w))=&\ 0\quad \forall\ w\in \mathcal{V},\\
  (\nabla \cdot u,q) - \nicefrac{1}{\lambda} (p,q)=&\ 0\quad \forall\ q\in \mathcal{U},\\
  (1- \kappa)(\varphi \sigma^+(u,p) : E_{\text{lin}}(u),\psi-\varphi) + G_c (-\nicefrac{1}{\epsilon} (1- \varphi),\psi-&\ \varphi)\\ 
  + G_c \epsilon (\nabla \varphi,\nabla(\psi-\varphi) )\geq&\ 0\quad \forall\ \psi \in \mathcal{K} \subset \mathcal{W},
  \end{aligned}
 \end{align*}
\end{problem}
where $\epsilon >0$ describes the bandwidth of the transition zone between broken and unbroken material.
This weak formulation in Problem~\ref{formMixed} can be reformulated to a complementarity system by introducing a Lagrange multiplier $\Lambda \in \mathcal{X}$, see~\cite{mang2019mesh,mang2019phase}.

The numerical treatment of the phase-field system in a monolithic fashion including the discretization as well as the adaptive refinement strategy are discussed in the following.


\subsection{Numerical Treatment and Programming Code}

Based on the complementarity formulation of Problem~\ref{formMixed}, with the help of a Lagrange multiplier, the crack irreversibility constraint is enforced, see~\cite[Section 4.1]{mang2019phase}. 
For the discretization in space, we employ a Galerkin finite element method in each loading step.
To this end, the domain $\Omega$ is partitioned into quadrilaterals. 
To fulfill a discrete inf-sup condition, Taylor-Hood elements with 
biquadratic shape functions ($Q_2$) for the displacement field $u$ and bilinear shape functions ($Q_1$)
for the pressure variable $p$ as well as for the phase-field variable are used.
For further details on the stable mixed form of the classical phase-field fracture model as well as the handling of the crack irreversibility condition and the numerical solving steps, we refer to~\cite{mang2019phase}.\\
The overall implementation is done in DOpElib~\cite{dope,DOpElib}
using the finite element library deal.II~\cite{dealII85}.


\subsection{Adaptive Refinement}

A residual-type a posteriori error estimator $\eta$ for the classical phase-field fracture model, presented and tested in~\cite{mang2019mesh}, provides a robust upper bound. Here, robust means that
the unknown constant in the bound does not depend on $\epsilon$ such
that the quality of the estimator is independent of $\epsilon$.
The mesh adaptation is realized using extracted local error indicators from the a posteriori error estimator in~\cite[Section 3.2]{mang2019mesh} on the given meshes over all loading steps.

In the following, $\mathfrak{M}^n$ denotes the mesh in the incremental step $n$ and $I^n_h$ is the corresponding nodal interpolation operator on the mesh $\mathfrak{M}^n$. 
The searched discrete quantities are denoted by an index $(\cdot)_h$, i.e., the displacement $u^n_h$, the phase-field variable $\varphi^n_h$, the pressure $p^n_h$, and the Lagrange multiplier $\Lambda^n_h$.
The adaptive solution strategy is given in the following.
\begin{Algorithm}\label{alg:refinement}
  Given a partition in time $t_0 < \ldots< t_N$, and an initial
  mesh $\mathfrak{M}^n = \mathfrak{M}$ for all $n = 0, \ldots, N$.
 \begin{enumerate}
 \item Set $\varphi^0_h = I_h^0 \varphi^0$ and
   solve the discrete complementarity system to 
   obtain the discrete solutions $\boldsymbol{u}^n_h,\varphi^n_h, p^n_h, \Lambda_h^n$ for all $n=1,\ldots,N$.
 \item Evaluate the error estimator in order to
   obtain $\eta^n$ for each incremental step.
 \item Stop, if $\sum_{n=1}^N(\eta^n)^2$ and $\|I_h^n
   \varphi^{n-1}-\varphi^{n-1}\|$ are small enough for all $n = 1,\ldots,
   N$.
\item For each $n=1,\ldots N$, mark elements in $\mathfrak{M}^n$ based on $\eta^n$ according to an
  optimization strategy, as implemented in deal.II~\cite{dealII85}.
\item Refine the meshes according to the marking and satisfaction of
  the constraints on hanging nodes.
\item Repeat from step 1.
\end{enumerate}
\end{Algorithm}


\section{Numerical Results}\label{Results}

In this section, the mixed phase-field model formulation is applied to simulate crack propagation in an L-shaped specimen with the help of adaptive refined meshes. 
First, the setup of the L-shaped panel test and the corresponding material and numerical parameters are given. Afterwards, the load-displacement curves and the crack paths are discussed for three different Poisson ratios from the standard
setting towards the incompressible limit $\nu \to 0.5$.

\subsection{Configuration of the L-shaped Panel Test}\label{config}
The L-shaped panel test was originally developed by Winkler~\cite{winkler2001traglastuntersuchungen} to test the crack pattern of concrete experimentally and numerically.
Concrete is compressible with a Poisson ratio of $\nu=0.18$. 
To simulate fracture propagation in nearly incompressible materials, within this work, the Poisson's ratio is increased towards an incompressible solid.

\begin{figure}[htbp!]
\begin{minipage}{0.47\textwidth}
\begin{tikzpicture}[xscale=0.65,yscale=0.65]
\draw (0,0) -- (2.5,0) -- (2.5,2.5) -- (5,2.5) -- (5,5) -- (0,5) -- (0,0);
\draw[white,fill=gray!30] (4,2.5) -- (5,2.5) -- (5,5) -- (4,5) -- (4,2.5);
\draw (2.5,0) -- (2.5,2.5);
\draw (2.5,2.5) -- (5,2.5);
\draw (5,2.5) -- (5,5);
\draw (5,5) -- (0,5);
\draw (0,5) -- (0,0);
\draw[<->] (0,5.5) -- (5,5.5);
\node at (2.5,5.73) {\scriptsize{$500 \mathrm{mm}$}};
\node at (2.5,5.2) {\scriptsize{$\Gamma_{\text{top}}$}};
\draw[<->] (-0.25,0) -- (-0.25,5);
\node at (-0.5,5.15) {\scriptsize{$y$}};
\node at (-1,2.5) {\scriptsize{$500 \mathrm{mm}$}};
\draw[<->] (0,-0.5) -- (2.5,-0.5);
\node at (1.25,-0.75) {\scriptsize{$250 \mathrm{mm}$}};
\node at (1.25,-1.2) {\scriptsize{$\Gamma_{\text{measured}}$}};
\node at (2.65,-0.5) {\scriptsize{$x$}};

\draw[<->] (5.25,2.5) -- (5.25,5);
\node at (6,3.75) {\scriptsize{$250 \mathrm{mm}$}};

\draw (0.1,0) -- (-0.05,-0.2);
\draw (0.3,0) -- (0.15,-0.2);
\draw (0.5,0) -- (0.35,-0.2);
\draw (0.7,0) -- (0.55,-0.2);
\draw (0.9,0) -- (0.75,-0.2);
\draw (1.1,0) -- (0.95,-0.2);
\draw (1.3,0) -- (1.15,-0.2);
\draw (1.5,0) -- (1.35,-0.2);
\draw (1.7,0) -- (1.55,-0.2);
\draw (1.9,0) -- (1.75,-0.2);
\draw (2.1,0) -- (1.95,-0.2);
\draw (2.3,0) -- (2.15,-0.2);
\draw (2.5,0) -- (2.35,-0.2);
\draw[blue,thick] (4.7,2.5) -- (5.0,2.5);
\draw[blue,->,thick] (4.75,1.95) -- (4.75,2.3);
\draw[blue,->,thick] (4.95,1.95) -- (4.95,2.3);
\node at (4.9,1.65) {\scriptsize{$\Gamma_{u_y}$}};
\node at (4.9,1.25) {\scriptsize{$30 \mathrm{mm}$}};
 \end{tikzpicture}
\captionof{figure}{Geometry and boundary conditions of the L-shaped panel test.}\label{lShaped}
 \end{minipage}
 \hfill
 \begin{minipage}{0.495\textwidth}
\renewcommand*{\arraystretch}{1.4}
\scriptsize
\begin{tabular}{|c|l|c|}\hline
\multicolumn{1}{|c}{Parameter} & \multicolumn{1}{|c}{Description}  & \multicolumn{1}{|c|}{Value} \\ \hline \hline
   $\mu$  & Lam\'{e} coefficient  & $10.95 \mathrm{kN/mm^2}$ \\ \hline
  $\lambda$  & Lam\'{e} coefficient & $6.16 \mathrm{kN/mm^2}$\\ \hline
   $\nu$ & Poisson's ratio & $0.18$ \\ \hline
   $G_c$ & Critical energy rate & $8.9\times 10^{-5} \mathrm{kN/mm}$ \\ \hline
    $h$  & Discretization parameter & $7.289 \mathrm{mm}$ \\ \hline
   $\epsilon$ & Bandwidth  & $14.0 \mathrm{mm}$ \\ \hline
   $\delta t$ & Incremental size &  $10^{-4} \mathrm{s}$ \\ \hline
      $\mathcal{I}$ & End time & $0.4 \mathrm{s}$\\ \hline
      $\kappa$ & Regularization parameter & $10^{-10}$ \\ \hline
 \end{tabular}
 \footnotesize
\captionof{table}{Standard settings of the material and numerical parameters for the L-shaped panel test.}\label{params}
\end{minipage}
\end{figure}

In Figure~\ref{lShaped}, the test geometry of the L-shaped panel test is declared.
In the right corner $\Gamma_{u_y}$ on a small stripe of $30 \si{mm}$ at the boundary, a special displacement condition is defined as a loading-dependent non-homogeneous Dirichlet condition:
 \begin{align*}
 u_y = t \cdot \si{mm/s},\quad \text{for}\ t \in I:= [0; 0.4 \si{s}],
\end{align*}
where $t$ denotes the total time and $T=0.4 \si{s}$ is the end time which corresponds to a displacement of $0.4 \si{mm}$. The time interval $I$ is divided into steps of the loading size $\delta t$.
\begin{remark}
To avoid developing unphysical cracks in the singularity on the
boundary $\Gamma_{u_y}$, the domain where the phase-field inequality
is solved, is constrained to the subset given by $x<=400 \mathrm{mm}$ similar to~\cite{mesgarnejad2015validation}. 
For $x>400 \mathrm{mm}$ we assume $\varphi=1$.
\end{remark}
In Table~\ref{params}, the required material and numerical parameters
for the L-shaped panel test are listed. Keep in mind, that the given
values for $\mu$ and $\lambda$ fit to the original material concrete
and are changed for other values of $\nu$ in
the following numerical tests, as listed in Table~\ref{poisson}. Further, the discretization parameter $h$ in Table~\ref{params} changes within the refinement steps, so $h$ is the starting mesh parameter on the coarsest mesh before adaptive refining.

\begin{figure}[htbp!]
\begin{minipage}{0.32\textwidth}
\renewcommand*{\arraystretch}{1.5}
\scriptsize
\begin{tabular}{|c|c|r|}\hline
\multicolumn{1}{|c|}{$\nu$} & \multicolumn{1}{|c|}{$\mu$} & \multicolumn{1}{|c|}{$\lambda$} \\ \hline \hline
$0.18$  & $10.95\cdot 10^3$ & $6.18\cdot 10^3$ \\ \hline
$0.40$ & $10.95\cdot 10^3$ & $42.36\cdot 10^3$ \\ \hline
$0.49$ & $10.95\cdot 10^3$ & $51.89\cdot 10^4$ \\ \hline
\end{tabular}
 \footnotesize
\captionof{table}{Tests with different Poisson's ratios.}\label{poisson}
\vspace{1cm}
\renewcommand*{\arraystretch}{1.4}
\scriptsize
\begin{tabular}{|l|c|r|}\hline
\multicolumn{1}{|l|}{$\nu$} & \multicolumn{1}{|c|}{min. $\#$DoF} & \multicolumn{1}{|c|}{max. $\#$DoF} \\ \hline \hline
$0.18$ & uniform & $213,445$\\ \hline
$0.18$ & $53,925$ & $125,599$ \\ \hline
$0.40$ & uniform & $213,445$\\ \hline
$0.40$ & $53,925$ & $121,709$ \\ \hline
$0.49$ & uniform & $213,445$\\ \hline
$0.49$ & $53,925$ & $91,574$ \\ \hline
 \end{tabular}
 \footnotesize
\captionof{table}{The minimal and maximal number of degrees of freedom (DoF) per incremental step on adaptive meshes.}\label{dofs}
\end{minipage}
\hspace{0.2cm}
\begin{minipage}{0.62\textwidth}
\begin{tikzpicture}[xscale=0.75,yscale=0.75]
\begin{axis}[
    ylabel = Load $F_y$ $\lbrack\mathrm{N}\rbrack$,
    xlabel = Displacement $u$ $\lbrack\mathrm{mm}\rbrack$,
 legend pos=north west, grid =major,
    x post scale = 1.1,
    y post scale = 1.43,
  ]
\addplot[]
table[x=Loading,y=018_uni_4,col sep=comma] {Data/lshaped_load_displ_uniform.csv};
\addlegendentry{$\nu = 0.18$ uniform}
\addplot[blue, densely dotted]
table[x=Loading,y=18_adaptive,col sep=comma] {Data/lshaped_load_displ_adaptive.csv}; 
\addlegendentry{$\nu = 0.18$ adaptive}
\addplot[cyan]
table[x=Loading,y=04_uni_4,col sep=comma] {Data/lshaped_load_displ_uniform.csv};
\addlegendentry{$\nu = 0.40$ uniform}
\addplot[green,dashed]
table[x=Loading,y=4_adaptive,col sep=comma] {Data/lshaped_load_displ_adaptive.csv};
\addlegendentry{$\nu = 0.40$ adaptive}
\addplot[magenta]  
table[x=Loading,y=049_uni_4,col sep=comma] {Data/lshaped_load_displ_uniform.csv};
\addlegendentry{$\nu = 0.49$ uniform}
\addplot[red,densely dotted]  
table[x=Loading,y=49_adaptive,col sep=comma] {Data/lshaped_load_displ_adaptive.csv};
\addlegendentry{$\nu = 0.49$ adaptive}
\end{axis}
\end{tikzpicture}
\caption{Load-displacement curves for the L-shaped panel test with different Poisson ratios and adaptively refined meshes versus uniform refinement. Weighted loading measured on the lower left boundary $\Gamma_{\text{measured}}$ labeled in Figure~\ref{lShaped}.}\label{loadDispl}
\end{minipage}
\end{figure}

\begin{figure}[htbp!]
 \centering
\begin{minipage}{0.3\textwidth}
 \includegraphics[width=0.95\textwidth]{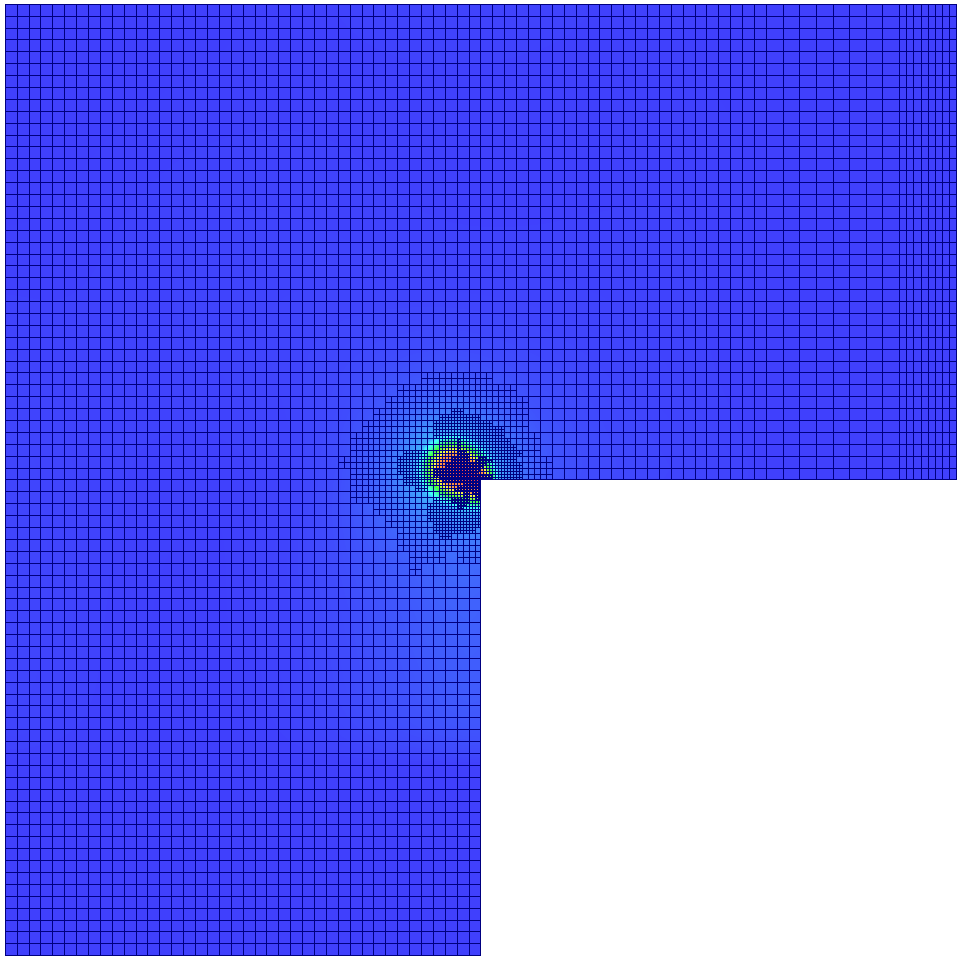}
\end{minipage}
 \hspace{0.1cm}
\begin{minipage}{0.3\textwidth}
 \includegraphics[width=0.95\textwidth]{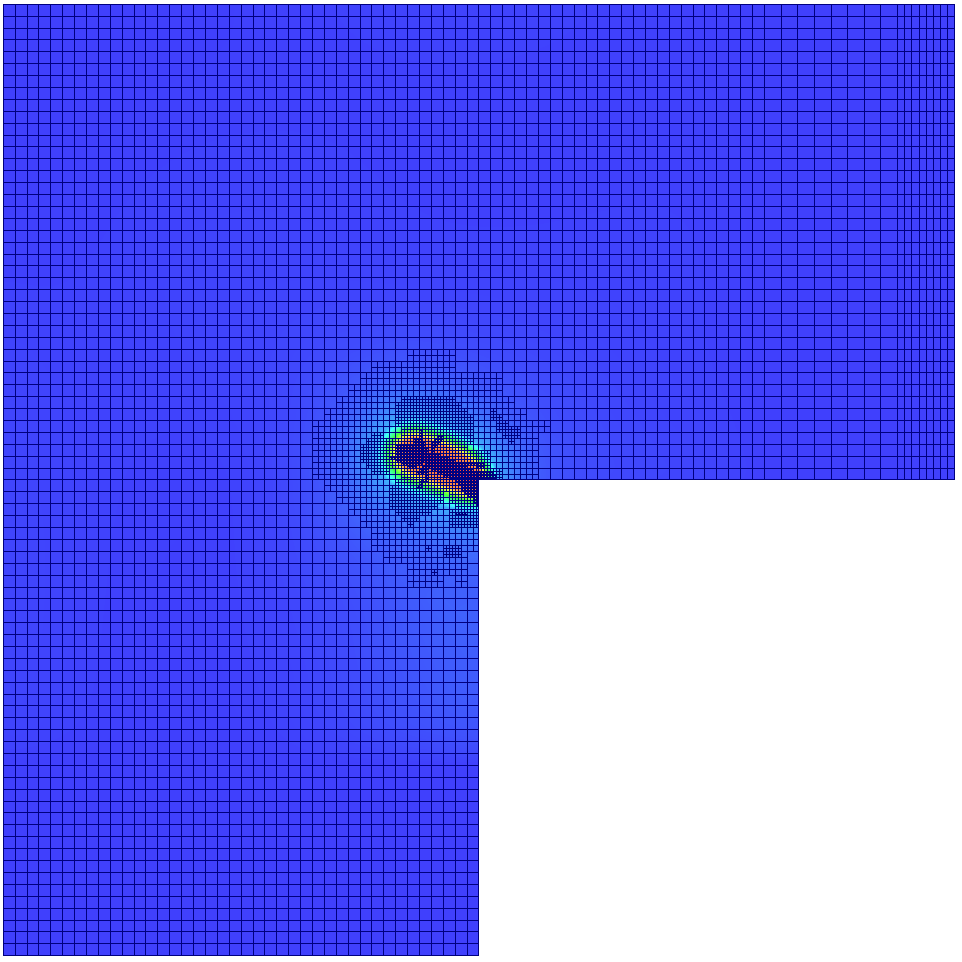}
 \end{minipage}
 \hspace{0.1cm}
\begin{minipage}{0.3\textwidth}
 \includegraphics[width=0.95\textwidth]{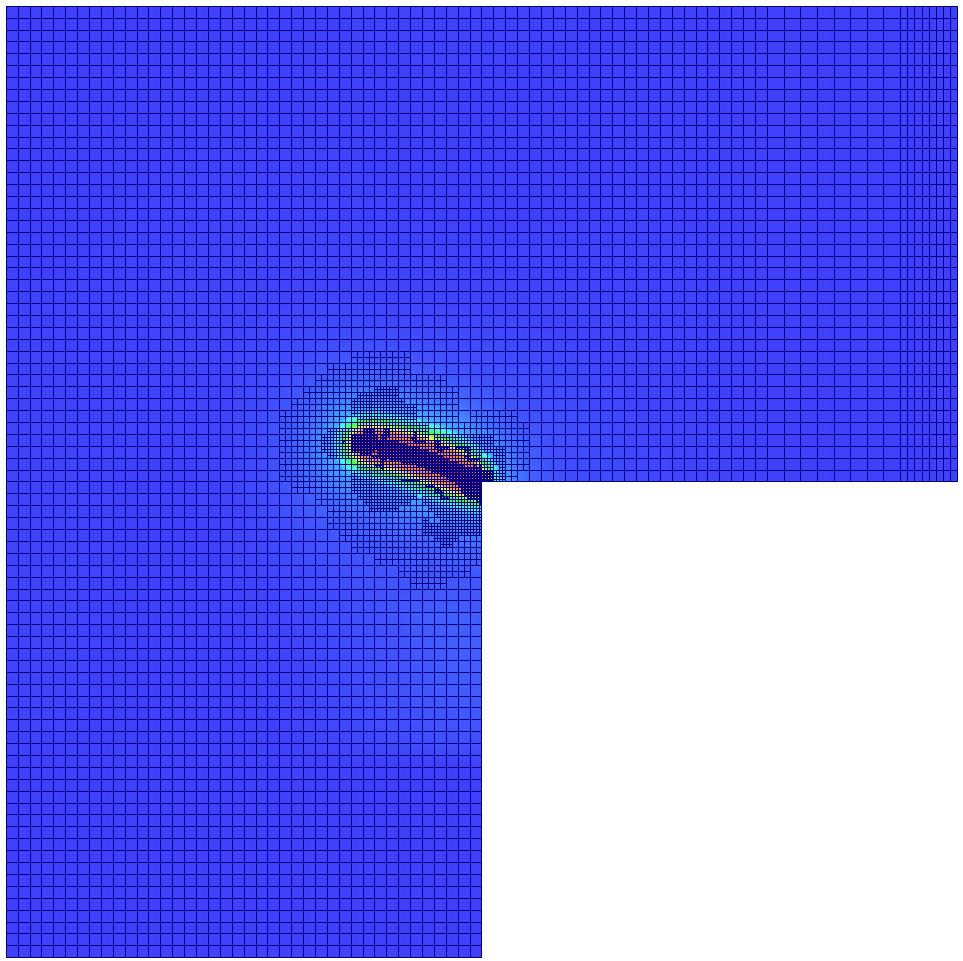}
 \end{minipage}\\
 \vspace{0.3cm}
\begin{minipage}{0.3\textwidth}
 \includegraphics[width=0.95\textwidth]{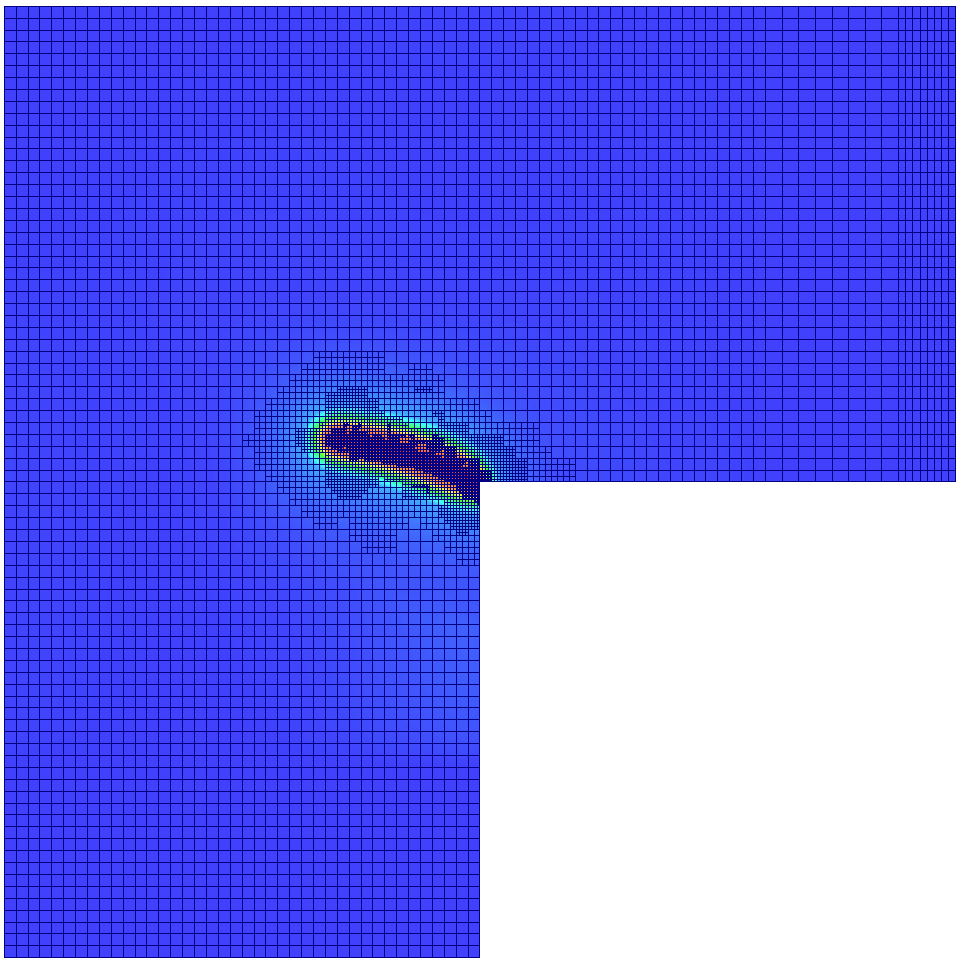}
\end{minipage}
 \hspace{0.1cm}
\begin{minipage}{0.3\textwidth}
 \includegraphics[width=0.95\textwidth]{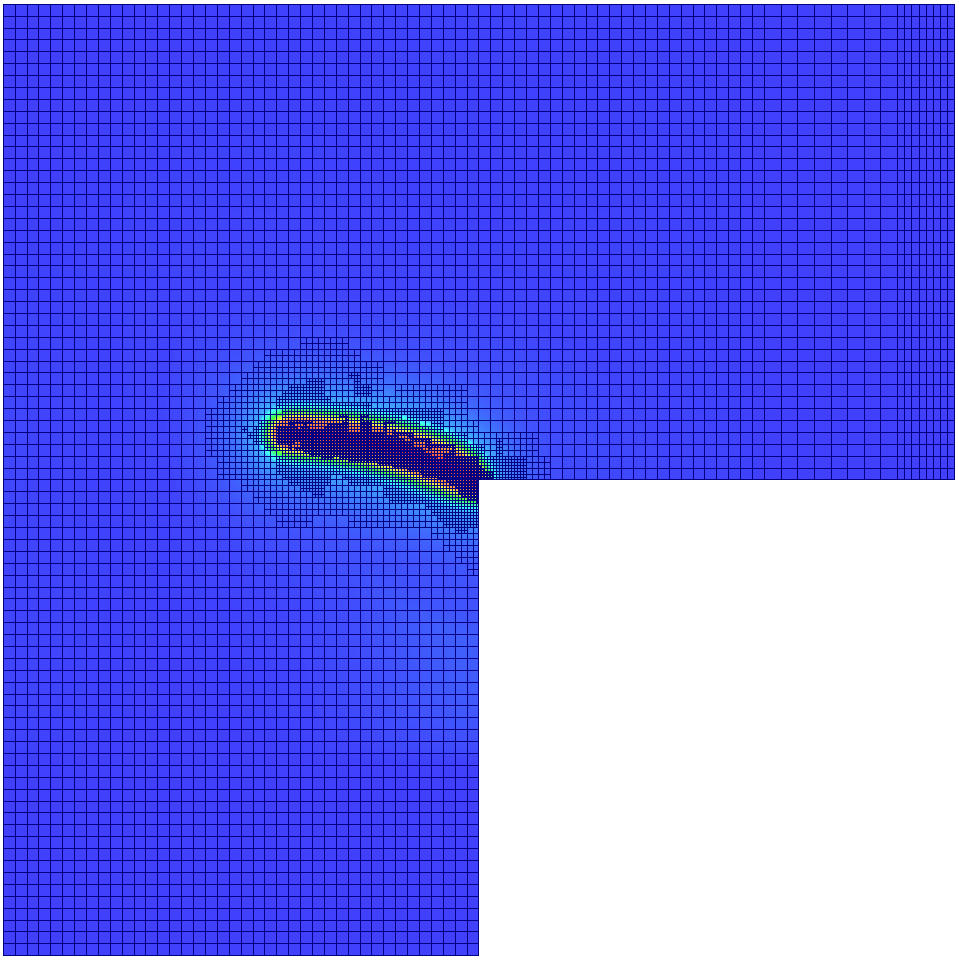}
 \end{minipage}
 \hspace{0.1cm}
\begin{minipage}{0.3\textwidth}
 \includegraphics[width=0.95\textwidth]{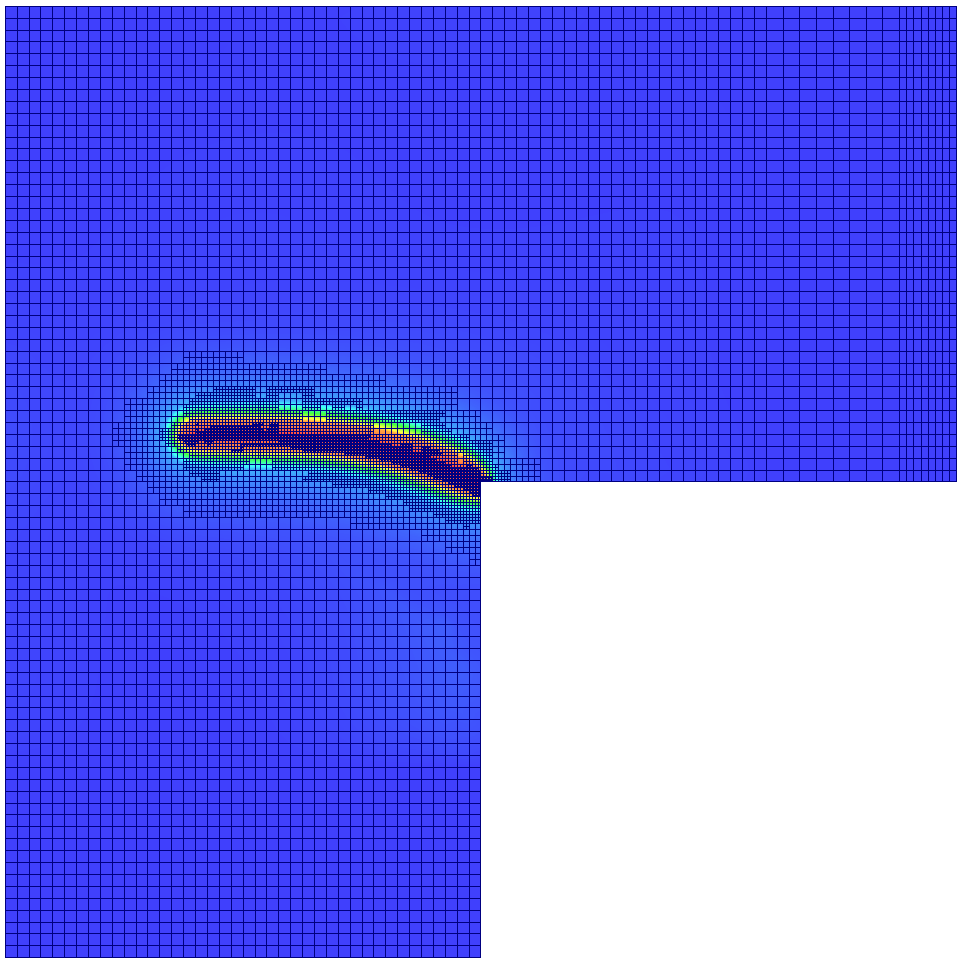}
 \end{minipage}
\caption{Poisson's ratio $\nu=0.40$. Snapshots of the phase-field function after three adaptive refinement steps in the incremental steps $0.2082, 0.209$, $0.2099$, $0.2136$, $0.2323$ and $0.2997 \mathrm{s}$ on the current adaptive mesh.}\label{plot3}
\end{figure}

\begin{figure}[htbp!]
 \centering
\begin{minipage}{0.29\textwidth}
 \includegraphics[width=0.75\textwidth]{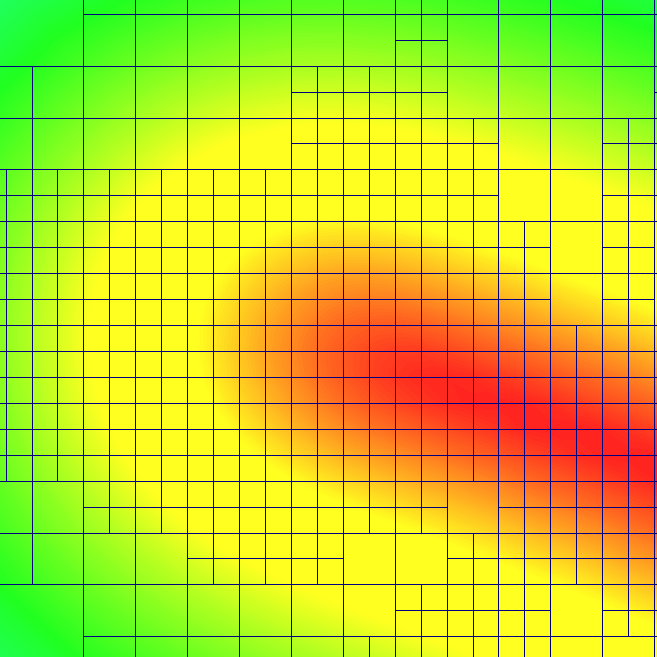}
\end{minipage}
 \hspace{0.4cm}
\begin{minipage}{0.29\textwidth}
 \includegraphics[width=0.75\textwidth]{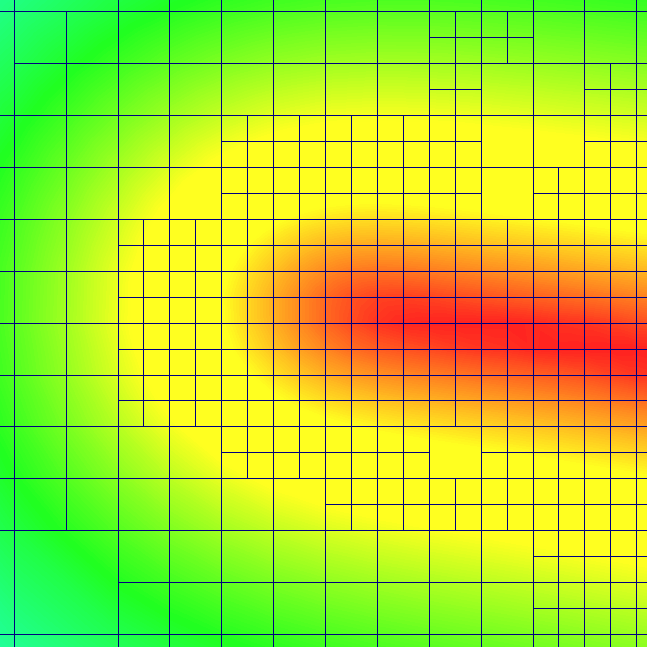}
 \end{minipage}
 \hspace{0.4cm}
\begin{minipage}{0.29\textwidth}
 \includegraphics[width=0.75\textwidth]{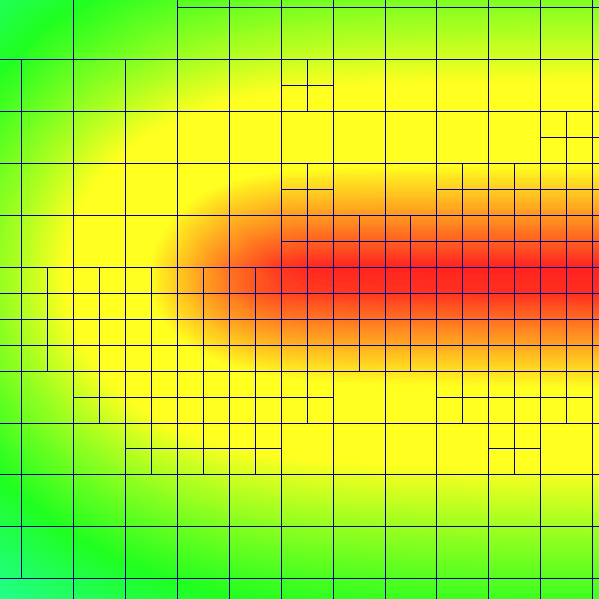}
 \end{minipage}
\caption{Poisson's ratio $\nu=0.40$. Enhanced extract of the phase-field function in the crack tip after three adaptive refinement steps in the incremental steps $0.2099$, $0.2176$ and $0.2997 \mathrm{s}$.}\label{plot4}
\end{figure}

In Table~\ref{dofs}, the minimal and maximal number of degrees of
freedom are given for three different test cases $\nu=0.18, \nu=0.40$
and $\nu=0.49$. The adaptive computations are based on a three times uniform refined mesh and three adaptive refinement steps.
For comparison, also the load-displacement curves for tests, executed on a four times uniform refined mesh, are added in Figure~\ref{loadDispl}.
The load-displacement curves in Figure~\ref{loadDispl} indicate that the higher the Poisson ratio, the higher is the maximal loading value before the crack starts propagating. 
Further, the path of the load-displacement curves for $\nu=0.18$, in particular for the adaptive test run in Figure~\ref{loadDispl}, coincide with the numerical and experimental results in concrete~\cite{winkler2001traglastuntersuchungen}.
In general, the adaptive computations exhibit a faster crack growth as it is expected in brittle materials, and may call for additional adaptive refinement of the time discretization for which models and indicators still need to be developed.
As a second quantity of interest, in Figure~\ref{plot3}, the crack path can be observed in certain incremental steps on adaptive refined meshes, exemplary for $\nu=0.40$.
The refinement strategy based on the error indicators steers the resolution of the crack area, especially of the crack tip as visualized in Figure~\ref{plot4}. 

\section{Conclusion}\label{Conclusion}
We have combined and extended~\cite{mang2019phase}
and~\cite{mang2019mesh} to adaptive refinement based on robust
residual-type a posteriori error estimators for phase-field model for fractures in nearly incompressible materials. 
The method is demonstrated on a numerical test for the L-shaped panel test.
Therefore, we proposed three test cases in Section~\ref{Results} with different Poisson ratios $\nu$ approximating the incompressible limit $\nu =0.5$. The load-
displacement curves of the three tests show a correlation between an increasing Poisson ratio and a stronger loading force.
In view of mesh adaptivity, we observed very convincing findings: the
mesh refinement is localized in the area of the (a priori unknown)
fracture path and allows to resolve the crack tip region.
Further, our adaptive refined meshes allow for a faster crack growth
compared to uniformly refined meshes.

\acknowledegment{This work was funded by the Deutsche Forschungsgemeinschaft
  (DFG, German Research Foundation) -- 392587580.
  It is part of the Priority Program 1748 (DFG SPP 1748) \textit{Reliable Simulation Techniques in Solid Mechanics. Development of
    Non-standard Discretization Methods, Mechanical and Mathematical Analysis}.

%

\bibliographystyle{abbrv}
\bibliography{lit}
\end{document}